\documentstyle[10pt,twoside]{siamltex}
\newtheorem{fed}[theorem]{Definition}
\newtheorem{teo}[theorem]{Theorem}
\newtheorem{lem}[theorem]{Lemma}
\newtheorem{cor}[theorem]{Corollary}
\newtheorem{pro}[theorem]{Proposition}
\newtheorem{rem}[theorem]{Remark}
\newtheorem{exa}[theorem]{Example}
\newtheorem{num}[theorem]{}
\def\dem{{\it Proof.\ }\rm}

\newfont{\bb}{msbm10}
\def\Bbb#1{\mbox{\bb #1}}

\def\zN{ \ \hbox{\rm N\hskip -10.8pt I
\hskip 4.7pt}}
\def\bm{\left(\begin{array}}
\def\em{\end{array}\right)}
\def\QED{\hskip 0.1truecm\vbox{\hrule\hbox{\vrule\hskip .2truecm
\vbox{\vskip .2truecm}\vrule}
\hrule}\hskip 0.1truecm}
\def\bull{\vrule height 1.0ex width .4ex depth -.1ex }
\def\inc{\subseteq}
\def\QED{\bull}
\def\<{\langle}
\def\>{\rangle}
\def\inv{^{-1}}
\def\la{\lambda}
\def\La{\Lambda}

\def\noi{\noindent}

\def\noi{\noindent}
\def\ha{Hadamard }

\def\ben{\begin{enumerate}}
\def\een{\end{enumerate}}
\def\IM{I(A)}
\def\IMJ{I(A_J)}
\def\IMm{I(A_m)}
\def\IMB{I(B)}
\def\IMBJ{I(B_J)}
\date{}

\begin{document}
\bibliographystyle{plain}

\title{ 
INDEX OF HADAMARD MULTIPLICATION BY POSITIVE MATRICES II {}}
\author {Gustavo  Corach \thanks
Supported UBACYT TX79, PIP 4463 (CONICET), Fundaci\'on
Antorchas and ANPCYT PICT 97-2259 (Argentina)
\and Demetrio Stojanoff \thanks Supported UBACYT TW49, 
PIP 4463 (CONICET) and ANPCYT PICT 97-2259 (Argentina)}

\vskip3truecm
 \maketitle 
\vskip1truecm
\noi
{\bf Gustavo Corach.}

\noi Instituto Argentino de Matem\'atica, 

\noi
Saavedra 15 Piso 3 (1083), Buenos Aires Argentina 

\noi
e-mail: gcorach@mate.dm.uba.ar

\bigskip

\noi
{\bf Demetrio Stojanoff}

\noi
Depto. de Matem\'atica, 

\noi
UNLP, 1 y 50 (1900), La Plata, Argentina  

\noi
e-mail: demetrio@mate.dm.uba.ar

\vglue2truecm

\begin{abstract}{ 
Given a definite nonnegative matrix $A \in M_n (\Bbb C )$, we
study the minimal index of $A$ : 
$ I(A) = \max \{\lambda \ge 0 :  A\circ B \ge \lambda B $ 
for all $ 0\le B\}$, where $  A\circ B$ denotes the Hadamard product
$( A\circ B)_{ij} = A_{ij} B_{ij}$. For any unitary invariant norm
$N$ in $M_n(\Bbb C )$, we consider the $N$-index of $A$:  
$I(N,A) = \min\{N(A\circ B) : B\ge 0$ and $N(B) = 1 \}$
If $A$ has nonnegative entries, 
then $I(A) =  I(\| \cdot \|_{sp}, A)$ if and only if there exists 
a vector $u$ with nonnegative entries such that 
$Au = (1, \dots, 1)^T$. We also show that
$I(\| \cdot \|_{2}, A)= I(\| \cdot \|_{sp}, {\bar A}\circ A)^{1/2}$. 
We give formulae for $I(N, A)$, for an arbitrary unitary invariant
norm $N$, when $A$ is a diagonal matrix or a rank 1 matrix.
As an application we find, for a bounded 
invertible selfadjoint operator $S$ on a Hilbert space, 
the best constant $M(S)$ such that $\|STS + S\inv T S\inv \| 
\ge M(S) \|T\| $ for all $0 \le T$.
}\end{abstract}

\vglue2truecm

\begin{keywords}
Hadamard product, definite nonnegative matrices, norm inequalities
\end{keywords}

\begin{AMS}
Primary 47A30, 47B15 \end{AMS}

\newpage
\section{Introduction}
We denote by $M_n = M_n (\Bbb C )$ the $C^*$-algebra of $n \times n$ matrices
over $\Bbb C $ and by $P(n) \subset M_n = \{  A \in M_n : A \ge 0 \}$ the
set of  definite nonnegative matrices. 
Given $A, B \in M_n$, we denote by $A \circ B$ their \ha
product, i.e. the matrix with entries $(A\circ B)_{i,j} = A_{ij} B_{ij}$.

For $A \in M_n$ we consider the map 
$\Phi_A: M_n \to M_n$ given by $\Phi_A(B)= A \circ B$,
for $B \in M_n$. By Schur Theorem \cite{Sch}, $A\circ B \in P(n) $ if $A, B \in P(n)$.
Thus, as a map between $*$-algebras, $\Phi_A $ is positive (actually it is
completely positive, \cite{pps} Prop. 1.2). 
Haagerup  \cite{uffe} has described the norm of $\Phi_A$ 
for $A \in M_n$ in terms of the factorizations $A = B^*C$. 
Other papers (see for example \cite{Ando}, 
\cite{AO}, \cite{cds}, \cite{cdp}, \cite{Da}, \cite{ho},
\cite{m}, \cite{pps} and  references included  therein) contain characterizations 
of several norms of $\Phi_A$. 
In the second author's paper \cite{St} 
the following problem was posed and partially solved: 
given $A \in P(n)$, calculate the infimum of $\| \Phi_A (B)\|$, for $B \in P(n)$
with norm one. This problem has two origins:
on one side it comes from the study of the index theory of 
completely positive maps on operator algebras (see \cite{BDH} and \cite{DH}).
On the other side, it was motivated by the search of optimal bounds for
certain operator inequalities (see \cite {CPR} and \cite{ACS}).

This paper is the natural continuation
of \cite{St} in both directions. 
We show several new characterizations of the two indexes
defined in \cite {St} and study the natural generalization of the index of
a positive matrix in terms of any unitary invariant norm in $M_n$.
We also get, for a  bounded selfadjoint 
invertible operator $S$ on   a Hilbert space ${\cal H}$, 
the best constant $M(S)$ such that 
$$
\|STS + S\inv T S\inv \| \ge M(S) \|T\| 
$$
for all positive operator $ T $ on ${\cal H}$. 
Let us give more explicit definitions of those \ha \ indexes: 

\medskip
\begin{fed} \rm \label{Imin}
The \ha {\bf minimal} index of $A \in P(n)$ is 
$$
\begin{array}{rl}
\IM      & = \max  \ \{ \ \lambda \ge 0 \ : \ A\circ B \ge \lambda B \quad 
\forall \ B \in P(n)  \ \} \\&\\
          & = \max \ \{ \ \lambda \ge 0 \ : \ \Phi_A - \lambda \ Id \ge 0  \quad
\hbox{ on } \quad P(n) \ \}  \\ &\\
          & = \max \ \{ \ \lambda \ge 0 \ : \ A-\lambda P \ge 0 \ \} 
\end{array} 
$$
where $P \in P(n)$ is the matrix with all its entries equal to 1. The
last equality follows from the fact that for $C\in M_n$, $\Phi_C \ge 0
\Leftrightarrow C \ge 0$ (see \cite{cdp}). 
In \cite{St} it is used the notation $II_A$ instead of $\IM$.

\end{fed}

\medskip
\begin{fed} \rm
For a unitary invariant norm $N$ in $M_n$, 
the \ha \ $N$-index for $A\in P(n)$ is
$$
\begin{array}{rl}
I(N,A)        & = \max \ \{ \ \lambda \ge 0 \ : \ 
N( A \circ B) \ge \lambda N(B) \quad \forall B \in P(n) \ \} \\
&\\
           & = \min \ \{ \ N( A \circ B ) \ : \ B \in P(n) \ \hbox{ and } \
N(B) = 1 \ \} \\ &\\
           & =  \min \ \{ \ N(B)^{-1} \ : \ 0 \ne B \in P(n) \ 
\hbox { and } \ N(A \circ  B ) \le 1 \ \} .
\end{array}
$$
For the Schatten p-norms, $1\le p < \infty$, 
we shall write $I(p, A) $ instead of $I(\| \cdot  \|_p , A)$. 
Note that the Schatten norm $\| \cdot \|_{\infty} $ is the spectral
norm $\|\cdot \|_{sp}$. The associated index will be denoted by 
$I(sp, \cdot )$. In \cite{St} it is used the notation $I_A$ instead of
$I(sp , A)$.
\end{fed}

\medskip
\begin{exa}\label{N1} \rm Let $A, B \in P(n)$. Then 
$$
\|B\|_{_1} = tr(B) = \sum_{i=1}^n B_{ii} \quad \hbox { and } \quad 
\|A\circ B\|_{_1} =tr(A\circ B)=  \sum_{i=1}^n A_{ii} B_{ii} .
$$
>From these identities it is easy to see that, for every $A \in P(n)$, 
$$
I(1, A) =  
\min _{1\le i \le n} A_{ii}.
$$
\end{exa}

\noi
We summarize the principal results of this work in the following theorems:

\smallskip

\begin{teo}
Let $A\in P(n)$. Then
\ben
\item $\IM \ne 0 $ if and only if the vector 
$p = (1, \dots , 1)^*$ belongs to the range of $A$. 
In this case, for any vector $y$
such that $Ay= p$, it holds 
$$\IM = \ \< Ay, y\> \inv = (\sum_{i=1}^n y_i )\inv = 
\min \ \{ \ \< Az,z\>  \ : \  \sum_{i=1}^n z_i = 1 \ \} .$$ 
(see  \ref{proIIA} and \ref {Bzz}).
\item Let $B \in P(m)$. Then $I(A\otimes B ) = \IM \ \IMB $ (see \ref{otim}). 
\item  If $A$ has nonnegative entries and $I(sp , A ) \ne 0$, then 
$$
I(sp , A ) = \IM 
$$ 
if and only if there exist a vector $u$ with {\bf nonnegative entries} 
such that $A(u) = p = ( 1, \dots , 1)^*$ (see Theorem \ref{i1=i2}).
\een 
\end{teo}
\smallskip
\begin{teo}
\ Let $A\in P(n)$. Then
\ben
\item $I(2, A) = I(sp ,  {\bar A}\circ A)^{1/2}$ (see Theorem \ref{elteo}).
\item If $A$ has nonnegative entries, then (see \ref{Bj})
$$ I(sp ,A ) = \min \{  I(sp ,A _J) :  J \subseteq \{1, \dots, n\}  \ \hbox{ and } 
 \ I(sp , A_J) = \IMJ  \}.$$ 
\een
\end{teo}

\smallskip
\begin{teo}
\ Let $A\in P(n)$ and $N$ an unitary invariant norm in 
$M_n$. Then 
\ben
\item $I(N, A) = 0 \Leftrightarrow $ $ A_{ii} = 0 $ for some $i = 1, 2, \dots , n$.
\item If $A $ has rank one, then
$I(N,A) = \min_{1\le i\le n} A_{ii}$ (see \ref{ran1}).
\item If $A$ is diagonal and invertible, then  
$I(N, A) = N'(A\inv) \inv , $ 
where $N'$ is the dual norm of $N$ (see \ref{lasN}). 
\een
\end{teo}
\smallskip

\begin{teo}
\  Let ${\cal H}$ be a Hilbert space and $S$ a bounded selfadjoint 
invertible operator 
on ${\cal H}$. Let $M(S)$ be the best constant such that 
$$
\|STS + S\inv T S\inv \| \ge M(S) \|T\| 
\quad \hbox{ for all }    \ 0\le T \in L({\cal H}). 
$$
Then  $M(S) = \min \{ M_1 (S) , M_2 (S)\}$, where
$$M_1(S) = \min_{\la \in \ \sigma (S)} \la ^2 + \la ^{-2}
\quad \quad and  
$$
$$ M_2(S)  =  \inf   \Big{\{} {(|\la | + |\mu |)^2 \over 1+\la^2 \mu^2} :
\la, \mu \in \sigma (S) , |\la | < |\mu | \hbox{ and }  
\la^2 \le {1\over |\la \mu |} \le \mu ^2 \Big{\}} . 
$$
In particular, 
if $\|S\|\le 1 $ (resp. $\|S\inv \|\le 1$), then 
$$M(S) = \|S\|^{2} + \|S \|^{-2} \quad  
(\hbox{resp. } \ \ \|S\inv \|^{2} + \|S\inv \|^{-2}).$$
\bigskip
\end{teo}

\section{The minimal index $\IM$}

Recall from Definition \ref{Imin} that 
for $A \in P(n)$,  
$$
\begin{array}{rl}
\IM      & = \max  \ \{ \ \lambda \ge 0 \ : \ A\circ B \ge \lambda B \quad 
\forall \ B \in P(n)  \ \} \\
          & = \max \ \{ \ \lambda \ge 0 \ : \ A-\lambda P \ge 0 \ \} 
\end{array} 
$$
where $P \in P(n)$ is the matrix with all its entries equal to 1. 
\bigskip
\noi 

\medskip
\begin{rem} \label{remIIA}
\rm Let $A \in P(n)$. 
\ben
\item The index $I(\cdot )$ is caalled  minimal because 
for every  unitary invariant norm $N$, it holds that
$\IM \le I(N, A)$. Indeed,  given  $B \in P(n)$, then
$A\circ B \ge \IM B $ and, by Weyl theorem, 
$\IM s_i (  B) \le s_i (A\circ B )$, $1\le i \le n$ ($s_i$ denote the 
singular values). Therefore $\IM N(B) \le N(A\circ B)$.
\item If $A$ is invertible, then (see Theorem 4.5 of \cite {St}) $\det (A+P) > \det (A)$ and 
$$
\IM = \frac { \det (A)}{\det (A+P) - \det (A)} = 
(\sum_{i,j=1}^n (A^{-1})_{ij})^{-1} = \ \< p, A\inv p\> \inv,
$$
where $p=(1, \dots, 1)^*$ and $P=pp^*$.
\item If a sequence $(A_m)_{m \in \zN}$ in $P(n)$ decreases  
to $A$, then (see Remark 4.3 of \cite{St})
$$
\lim_{m \to \infty} \IMm = \inf _{m \in  \zN } \IMm= \IM .
$$
\item If $J \inc \{ 1, \ 2, \dots, \ n\}$ and $A_J$ is the principal submatrix of 
$A$ associated to $J$, then $\IM \le \IMJ$. 
Also $I(N, A)\le I(N, A_J)$ for every unitary invariant norm $N$. 
Indeed, these inequalities can be deduced easily from the definitions of the 
index.
\een
\end{rem}
\bigskip

\noi
By the Remark above, if $A \in P(n)$ is invertible and $y = A\inv p$, then  
$$
0\ne \IM =\ \< p, A\inv p \> \inv \ = \ \< Ay, y\> \inv .
$$
In the next Proposition we shall see that the same formula also holds for
non invertible $A\in P(n)$, under the hypothesis that $p$ belongs to the
range of $A$.
\medskip
\begin{pro}\label{proIIA}
Let $A \in P(n)$. Then  $\IM \ne 0 $ if and only if the vector 
$p = (1, \dots , 1)^*$ belongs to the range of $A$. 
In this case, for any vector $y$
such that $A(y)= p$, we have that 
$$
\IM = \ \< Ay, y\> \inv = (\sum_{i=1}^n y_i )\inv \ . 
$$
\end{pro}
\noi \dem
First note that if 
$p$ lies in the range of $A$ then  $p \in \ker A ^\perp$. This means that
$\< Aw,w\>  \ = \ \< Ay,y \> $ for every pair $w, y$ such that $Aw = Ay= p$. 

Let $Q$ the orthogonal projection onto $\ker A$. 
Then the sequence $A_m = A+\frac1m Q$, $m \in \zN$, decreases to $A$   
and  $\IM = \ \lim_{m\to \infty}\IMm$. 
Note that $A_m$ is invertible for all $m \in\zN$.

If there exists a vector $y$ with $Ay = p$, let $y = w+z$ with 
$z \in \ker A$ and $w \in \ker A ^\perp$. 
Note that  $A_m\inv \ p= w , \ \forall \ m$. Therefore 
$$
\begin{array}{rl} 
\IM &= \ \lim_{m\to \infty}\IMm \\ & \\
& = \  \lim_{m\to \infty} \< p, A_m\inv p\> \inv \\ & \\
& =  \ \< p,w\>  \inv\ =   \ \< p,y\> \inv  \\ & \\ 
&  = \ \< Ay,y\>  \inv \ne 0. \end{array}
$$

\noi 
On the other hand, if $p \notin Im A$, let 
$$p = y+z \ , \quad y \in Im A  \quad \mbox{  and } \quad 0\ne z \in \ker A.
$$
If $Ax = y$ with $x \in \ker A^\perp$, then 
$A_m\inv p =x + m z$. Therefore 
$$
\IMm\inv = \ \< p,A_m\inv p \> \ = 
\ \< p,x\>  + m\< p,z\>  \ =\ \< p,x\>  + m\|z\|^2 \to \infty .
$$
Then $\IM = 0$ \QED

\medskip

\begin{cor}\label{otim}
Let $A \in P(n)$ and $B \in P(m)$. Consider the Kronecker product
matrix 
$$A\otimes B = \left(\begin{array}{ccc}
A_{11}  B & \dots & A_{1n} B \\
\vdots    &   \vdots & \vdots \\          
A_{n1}  B & \dots & A_{nn} B 
\end{array}\right)  \in P(nm).
$$ 
Then
$$
I(A\otimes B ) = \IM \ \IMB .
$$
\end{cor}
\noi \dem Suppose that $\IM \ne 0 \ne \IMB$. For any $k \in \zN$ denote by 
$p_k \in \Bbb C ^k$ the vector with all its entries equal to $1$.   
Let $x\in \Bbb C ^n $ and $y \in \Bbb C  ^m$ such that $Ax = p_n$ and $By = p_m$.
Then $z = (x_i y , x_2 y , \dots , x_n y )^t \in \Bbb C ^{nm}$ verifies that 
$(A\otimes B )z = p_{nm}$. Also 
$$
\begin{array}{rl}
I(A\otimes B )\inv  & = \ \< z, p_{nm}\>  \\&\\ 
        & = \sum_{i=1}^n (x_i \sum _{j=1}^m y_j) \\&\\
		& = \ \<  x, p_n  \>  \ \< y,  p_m \>  \ = 
		\IM \inv \ \IMB\inv  . \end{array}
$$
If $\IM =0 $ or $\IMB = 0$ then $I(A\otimes B ) = 0$: this can be verified 
by just multiplying by appropiate matrices of the type $C\otimes D$, since 
$(A\otimes B) \circ (C\otimes D) = (A\circ C) \otimes (B\circ D) $\ \ \QED

\medskip
\begin{rem} \rm 
As a particular case of Corollary \ref{otim}, the inflation
matrix $A^{(m)} = p_m p_m^* \otimes A$ verifies $I(A^{(m)}) = \IM $ for all
$A\in P(n)$ and $m \in \zN$. Using proposition (3.9) of \cite {St}, it can be also 
shown that 
\begin{equation}\label{infl}
I(sp , A^{(m)}) = I(sp , A) 
\end{equation}
for all
$A\in P(n)$ and $m \in \zN$. Indeed, if $A = BB^*$, then
$$
I(sp , A) = \min_{\|x\|=1}   \|D_x B\|^2  , 
$$
where $D_x$ denotes the diagonal matrix with the vector $x$ in its diagonal.
Using also this formula for $ A^{(m)}$ and the fact that 
$$ 
A^{(m)} = \bm {cccc} B& 0 &\dots &0\\
B& 0 &\dots &0 \\
\vdots & \vdots& \vdots&\vdots \\ B& 0 &\dots &0 \em  \ 
\bm {cccc} B^* & B^*  &\dots &B^* \\
0& 0 &\dots &0 \\
\vdots & \vdots& \vdots&\vdots \\ 0& 0 &\dots &0 \em   , 
$$
one easily gets the equality (\ref{infl}).
\end{rem}

\medskip
\begin{cor}\label{Bzz}
Let $A\in P(n)$. Then
$$
\IM = \min \ \{ \ \< Az,z\>  \ : \  \sum_{i=1}^n z_i = 1 \ \}
$$
\end{cor}
\noi \dem 
Given a vector $z$ such that  $\< z,p\>  \ = 1$, then
$$
\begin{array}{rl}
\< A z,  z\>   \ = \ \sum_{ij}A_{ij} z_j{\bar z_i}
   &  = \ \< (A\circ {\bar z}{\bar z}^*) p, p \>  \\&\\ 
   & \ge \ \IM \< {\bar z}{\bar z}^* p, p\>  \ 
      = \ \IM  \sum_{ij}{\bar z_i}{\bar z_j}  \\&\\ 
   &  = \ \IM \  \< p, z\> ^2 \ = \IM \end{array}
 $$
If there exists $ x \in \Bbb C  ^n$ such that $Ax = p$, the vector $z = \IM x $ 
verifies 
$$
\< p, z\>  \ = \IM \ \< p, x \>  \ = 1 \quad  \hbox { and }  
$$
$$
\< Az, z \>  \ = \ \IM  \< p, z\>  \ = \IM 
$$
by Proposition \ref{proIIA}. But if $p \notin (\ker A)^\perp $, then
there exists $z \in \ker A$ such that $\< z, p \>  \ = 1$ and $\< Az,z\>  \ = 0 = \IM$ \QED

\bigskip
The following result seems to be very well known. We shall state it with 
a proof for the sake of completeness and because we shall use it in a 
precise formulation (in Theorem \ref{i1=i2} and Proposition \ref{el x}).

\medskip
\begin{lem}\label{el z} Let  $ B \in P(n) \cap M_n (\Bbb R )$ and 
$M = \{ z \in \Bbb R^n : \sum_i z_i = 1 \ \}$. Then
$$
V_1 = \{ z \in M  :  \< Bz, z \>  \ = \IMB \} = 
\{  z \in M  :  Bz = \IMB  p  \} =  V_2  \ne \emptyset ,
$$
where $p = ( 1, \dots , 1 )^* $. Moreover, any local extreme point of the map 
$G(z) = \ \< Bz,z\>  $ restricted to $M$ belongs to $V_2 $.
\end{lem}

\noi \dem 
It is clear that $V_2 \subseteq V_1$. 
Recall from the proof of Corollary \ref{Bzz}
(and the fact that $B \in M_n(\Bbb R )$) 
that $V_2 \ne \emptyset$, so $\IMB = \min \{ \ \< Av,v\>  \ : \  v \in M \ \}$. 
The map $G(z) = \  \< Bz, z\>  \ = \sum_{i,j}b_{ij} z_jz_i$ is differentiable.
So we can use Lagrange multipliers in order to find its critical points in $M$.  
Let $F(z, \la) =  \sum_{i,j}b_{ij} z_jz_i - 2 \la (\sum _1^n z_i - 1)$. 
Then, since $B^t = B$,  
$$
{\partial F \over \partial z_i } (z, \la) = 
\sum _{j=1}^n b_{ij}z_j  + \sum _{j=1}^n b_{ji}z_j -2 \la  
=2 \sum _{j=1}^n b_{ij}z_j  -2\la  
= 2 ((Bz)_i - \la ).
$$
So, if $z\in M$, ${\partial F \over \partial z_i } (z, \la) = 0 $ for all $i$ if and only if 
$Bz = \la p$. 
But, in that case,  
$$ \IMB \le \ \< B z, z\>  \ = \la \< p,z\>  \ = \la.$$ 
If $\IMB = 0$ then $\la = 0 $,  since $p \notin $ Im $B$, by 
Proposition \ref{proIIA}.
If $\IMB > 0$ then also $\la = \IMB$, since $y = \la \inv z $ verifies
$By = p$ and 
$$
\la = \ \< Bz, z\>  \ = \la^2  \< By, y\>  \ = \la ^2 \IMB \inv.
$$
So $z\in M$ is a critical point if and only if \ $z \in V_2 $. Since each
local extreme mast be a critical point, this shows that $V_1 \subseteq V_2$ 
and the final assertion. \ \ \QED

\bigskip
The following Lemma, which is rather elementary, is useful
in order to identify vectors $x$ such that
some index is attained at the matrix $xx^*$.

\medskip
\begin{lem}\label{autovalor}
Let $A \in M_n$ and  $x \in \Bbb C  ^n$ with $\|x\|=1$. 
Let $y = x \circ {\bar x} = (|x_1|^2, \dots , |x_n|^2)^*$. 
Denote by $p = (1, \dots , 1)^*$. Then 
\ben
\item  If $Ay = \lambda p$, with $\la\in \Bbb C $, then 
$x$ is an eigenvector of the matrix $A  \circ xx^*$. 
\item If all $x_i \ne 0$ 
and  $(A  \circ xx^* )x = \la x $ for some $\la\in \Bbb C $, then 
$Ay= \la p$. 
\een
\noi If $A \in P(n)$, the eigenvalue $\la$ associated 
to $x$ must be $\IM$ and 
$Ay = \IM \ p$.
\end{lem}
\dem Suppose that $Ay = \lambda p$. Then 
$$
\begin{array}{rl}
(A\circ xx^* )x 
& = (a_{ij}x_i \bar{x}_j) 
\left( \begin{array}{c} x_1 \\ \vdots \\x_n  \end{array} \right) 
= \left( \begin{array}{c} ( \sum_j a_{1j} |x_j| ^2)  \ x_1 \\ \vdots \\
(\sum_j a_{nj} |x_j| ^2 ) \ x_n \end{array} \right) \\&\\
& = \left( \begin{array}{c}  (Ay)_1  \ x_1 \\ \vdots \\
(Ay)_n  \ x_n \end{array} \right) = \lambda \ x  .\end{array}
$$
The same equation shows that if $x$ is  an eigenvector of $A  \circ xx^*$
with non zero entries, then  $Ay =\la p$, where $\la$ is the 
eigenvalue of $x$.
If  $A\in P(n)$ and $\IM = 0$, then $\lambda = 0$ since $p \notin $ Im$A$.
If $\IM \ne 0$, then 
$p \in $ Im $ A = (\ker A)^\perp$. So $Ay \ne 0 $ since $1 = \|x\|^2  
= \ \< p, y \>  \ \ne 0$. Then  
$\la \ne 0$. If $z = \la \inv y$, then $Az = p$ and 
$$
1 = \ \< p,y\>  \ \ = \la  \ \< Az,z\>    \ \ =  \ \la \ \IM \inv ,
$$
by Proposition \ref{proIIA}. So $\lambda = \IM$ \QED

\bigskip

\noi Concerning the problem of characterize those matrices $A\in P(n)$ such that $\IM = 
I(sp , A)$, in \cite{St} it is shown that
for $A = \left( \begin{array}{cc} a&b\\
\bar b& c \end{array} \right) \in P(2)$, then 
\begin{equation}\label{bmin}
0 \ne I(sp, A) = \IM \quad \Leftrightarrow  \quad b \in \Bbb R \ \hbox { and } \ 
 0 \le b  \le  \min \{a, c\} \ne 0. 
\end{equation}
This is easily seen to be equivalent to the conditions
\ben
\item $A$ has nonnegative entries.
\item There exists a vector $z$ with nonnegative entries
such that $Az=(1,1)^*$ (if $A$ is invertible, this means that $A\inv (1,1)^*$
has nonnegative entries).
.
\een

\noi In the next Theorem we prove that, for  matrices 
of any size with nonnegative entries, condition 2 is equivalent to 
the identity $I(sp, A)= \IM$. 

\bigskip
\begin{teo}\label{i1=i2} 
Let $A\in P(n)$ with nonnegative entries such that all $A_{ii} \ne 0$. 
Then the following conditions are equivalent:
\ben
\item There exist a vector $u$ with {\bf nonnegative entries} 
such that $Au = p = ( 1, \dots , 1)^*.$
\item $I(sp , A ) = \IM $.
\een
In that case, if \ $y = \IM \ u $, then
\ben
\item [(a)]Let $x = (y_1^{1/2}, \dots , y_n ^{1/2})^t $. Then 
$\|x\|=1$ and  
$$
\| A\circ xx^*\| = I(sp , A) 
$$
\item [(b)]Let $J = \{ i : u_i \ne 0 \}$ and  denote by $A_J$ the 
principal submatrix of $A$
determined by $J$. Then  $\IM = \IMJ = I(sp ,A_J) = I(sp , A)$.
\een
\end{teo}

\dem 
\noi {\bf 2 $\to$ 1} 
Suppose $\IM = I(sp , A) $. Note that $I(sp, A) \ne 0$ by Corollary 3.4 
of \cite{St} (or Remark \ref{Rem1}).
Let $M = \{ z \in \Bbb R^n : \sum_i z_i = 1 \ \}$.
By Lemma 2.1 of \cite{St}, there exists 
$x\in \Bbb C ^n$ such that  $\|x\| = 1$ and $I(sp, A)= \|A \circ xx^*\|$. 
Denote $y = (|x_1|^2, \dots , |x_n|^2)^* \in M$. Then, by Corollary \ref{Bzz}, 
$$ 
\begin{array}{rl}
\IM = I(sp, A) & \ge \ \<  (A \circ xx^*)x, x\>  \\&\\ 
               & = \sum _{ij} A_{ij} |x_i|^2 |x_j|^2 \\&\\
               & = \ \< Ay, y\>  \ \ge \IM . \end{array}
$$
By Lemma \ref{el z}, $Ay = \IM p$ and $y$ 
has nonnegative entries. Take $u = \IM \inv y$.

\noi {\bf 1 $\to $ 2} Let $u$ a vector with nonnegative entries such that 
$Au = p$. Let  $y = \IM u \in M $ and $x$ as in item (a). 
Clearly  $ \|x\| = 1 $.  
By Lemma \ref{autovalor} we know that $x$ is an eigenvector of 
$A\circ xx^*$ with eigenvalue $\IM$. Recall that always $\IM \le I(sp ,A)$.

\noi{\bf Case 1.} Suppose that $x$ has strictly positive entries.  Since 
$A\circ xx^*$ has  nonnegative entries, it is well known (see 
Corollary 8.1.30 of \cite{HJ}) that the eigenvalue $\IM$ of $x$ must be the
spectral radius of $A\circ xx^*$. Since $A\circ xx^*\in P(n)$  we deduce that
$\IM = \|A\circ xx^*\| \ge I(sp ,A)$. 

\noi {\bf Case 2.} Let $J = \{ i : x_i \ne 0 \}$,  $A_J$ the principal 
submatrix of $A$ determined by the indexes of $J$ and similarly define $x_J$. 
Then $x_J$ is an eigenvector of $A_J \circ  x_Jx_J^*$ with eigenvalue $\IM$.
Note also that $A_J \circ  x_Jx_J^* \ge \IMJ  x_Jx_J^* $ and 
$x_Jx_J^* (x_J)   = \|x_J \|^2  x_J = x_J$.  Then 
$$
0 \le \< (A_J \circ  x_Jx_J^* - \IMJ  x_Jx_J^*) x_J, x_J\>  \ = \IM- \IMJ 
$$
and, by Remark \ref{remIIA}, $\IMJ = \IM$. Now, as in case 1, we can deduce that 
$$
\IM = \IMJ = \|A_J\circ x_Jx_J^*\| \ge 
I(sp ,A_J) \ge I(sp , A),
$$ 
where the last inequality holds by Remark \ref{remIIA}. 
Clearly  $I(sp , A)$ is  attained at $xx^*$ \ \ \QED

\medskip
\begin{rem} \rm In the last Theorem, the hypothesis that all $A_{ij} \ge 0$ is
essential in the implication 1 $\to$ 2. Indeed, consider $A = \bm {cc} 2&-1\\-1&1 
\em $ and $u = (2, 3)^*$. Then $Au = (1,1)^*$ but $1/5 = \IM \ne I(sp, A) = 
1$. On the other hand 2 $\to $ 1 remains valid without the mentioned hypothesis (our proof
only uses that $A\in M_n(\Bbb R )$ and the general case 
follows by just extending Lemma \ref{el z} to the complex case).
In any case we conjecture that condition $\IM = I(sp, A)$ actually 
implies that all $A_{ij} \ge 0$, as in the $2 \times 2$ case.
\end{rem}

\section{$I(sp, A)$ and $I(2, A)$}

In this section we shall study the relation between the \ha \ 
indexes associated to the spectral and the Frobenius norms. 

In  Lemma 2.1 of \cite {St} it is shown that the index $I(sp , \cdot) $ 
is always attained at rank one projections. It is natural to 
conjecture that the same result holds for any unitary invariant norm $N$.
In the following Proposition we show 
that our conjecture is true for the Frobenius norm:   
\medskip
\begin{pro}\label{uno}
Let $A\in P(n)$. Then there exists $x \in \Bbb C ^n$, 
$\|x\|=1$ such that $I(2,A) =\|A \circ xx^*\|_2$. That is, $I(2,A) $ 
is attained at rank one projections.
\end{pro}
\smallskip
\dem 
Let $\lambda = \max \{\mu \ge 0 :\| A \circ B\|_2 \ge \mu \|B\|_2 \ 
\forall B \in P(n) $ of rank one $\}$. 
By its definition $\lambda \ge I(2,A)$. Let us prove that
$\|A\circ B\|_2 \ge \lambda \| B\|_2 $ for all $B \in P(n)$. Indeed, 
for $ B \ge 0$,  write $B= \sum_{i=1}^k B_i$ where each $B_i$ has 
rank one, $B_i \in P(n)$  and $B_i B_j =0 $ if $i\neq j$. So 
$$
\lambda ^2 \|B\|_{_2}^2 = \lambda ^2 \sum_{i=1}^k 
\|B_i\|_{_2}^2 \le  \sum_{i=1}^k \|A\circ B_i\|_{_2}^2 .
$$
The proof is completed by showing that 
$$\sum_{i=1}^k \|A\circ B_i\|_{_2}^2 \le \|A\circ B\|_{_2}^2 .
$$
Indeed, by the parallelogram law (actually its generalization to $k$ 
vectors in the Euclidean space $(M_n, \| \cdot \|_2 )$ ), 
\begin{equation}\label{paral}
\sum_{\sigma} \|\sum_{i=1}^k \ (-1)^ {\sigma_i} \ A\circ B_i\|_{_2}^2 = 
2^k \sum_{i=1}^k \|A\circ B_i\|_{_2}^2 ,
\end{equation}
where the sum runs over all $\sigma \in \{0, 1\}^k$.
But clearly 
$$
-A\circ B \le \sum_{i=1}^k \  (-1)^ {\sigma_i} \ A\circ B_i \le A\circ B 
\ , \quad \hbox{ so }
$$
$$
|\<   \Big( \sum_{i=1}^k \ (-1)^ {\sigma_i} \ A\circ B_i \Big) \ x , x\> | \ \le 
\ \  \< (A\circ B) x, x\>  \ , \quad x\in \Bbb C ^n
$$
for all such $\sigma$. Therefore 
$\|\sum_1^k \ (-1)^ {\sigma_i} \ A\circ B_i \|_{_2}^2 \le \|A\circ B\|_{_2}^2 $
for all $\sigma$. Then looking at the convex combination which 
follows from equation (\ref{paral}), we conclude that
$\sum_1^k \|A\circ B_i\|_{_2}^2 \le \|A\circ B\|_{_2}^2$ and the proof
is complete \QED

\medskip
\begin{pro}\label{el x}
Let $A \in P(n)$. Then 
\ben
\item There exists a vector $x$ with nonnegative entries such that 
$\| x\|=1$ and $\| A\circ xx^*\|_2 = I(2,A)$.
\item $x$ is an eigenvector of the matrix $A \circ \bar{A} \circ xx^*$ with 
eigenvalue $\IMBJ$, where $B=A\circ {\bar A}$ and  $J = \{ i : x_i \ne 0 \}$.
\een
\end{pro}
\dem Let $y$ be an unit vector such that $\|A\circ yy^*\|_2 = I(2,A)$. 
Let $y_i ={ \bar w}_i |y_i|$, $|w_i|=1$, $1 \le i \le n$. If $w =
(w_1, \dots , w_n)^t$ and $D_w$ is the
diagonal and unitary matrix with the vector $w$ in its diagonal, then for
each $C \in M_n$, 
\begin{equation}\label{mod}
ww^* \circ C = (w_i \bar {w}_j \ C_{ij}) = D_w \ C \ D_w^*. 
\end{equation}
Therefore the \ha product by $ww^*$ doesn't change 2-norms. 
Denote by $x = (|y_1|, \dots , |y_n|)^*$. Then
$$
I(2,A) = \|A\circ yy^*\|_2 = \|A\circ yy^*\circ ww^*\|_2  =\|A\circ x x^*\|_2 
$$
showing item 1. Let $B= A \circ \bar{A} \in P(n)$.
Let $y \in \Bbb R_+ ^n$ with $\|y\|=1$ and let $z=(y_1^2, \dots , y_n^2)^*$. Then 
$$ \|A\circ yy^*\|_{_2}^2 = \sum_{i,j} |a_{ij}|^2 y_i^2 y_j^2
= \sum_{i,j}b_{ij} z_jz_i = \ \< Bz, z\>  
$$
and $\sum_1^n z_i = 1$. Then $\|A\circ yy^*\|_{_2} = I(2, A) $ if and only if 
$\< Bz, z\> $ is the minimum of the map $G(v) = \ \< B v, v\>  $ in the simplex
$\Delta = \{ v \in (\Bbb R_{\ge 0})^n :  \ \sum _1^n v_i = 1 \ \}$.
Using Lemma \ref{el z}, we know that if $z $ belongs to the  interior 
$\Delta ^\circ $ of $\Delta$, then $z $ is a local extreme 
of $G$ in the plane $M = \{ z \in \Bbb R^n : \sum_i z_i = 1 \ \}$, so 
$Bz = \IMB \ p$.

\noi
If the vector $x$ of item 1 verifies that $x_i > 0 $ for all $i$, then 
$z = x \circ x \in \Delta ^\circ $ and $Bz = \IMB \ p$. 
Then item 2 follows from Lemma
\ref{autovalor} with eigenvalue $\IMB$.
If some $x_i = 0$, let $J = \{ i : x_i \ne 0 \}$,  $B_J$ the 
principal submatrix of $B$
determined by the indexes of $J$ and similarly define $x_J$. Then 
$I(2, A) = \|A \circ xx^*\|_2 = \|A_J \circ x_Jx_J^*\|_2 \ge I(2, A_J)$
and 
$$
I(2,A ) = I(2,A_J )=\|A_J \circ x_Jx_J^*\|_2,
$$
since the other inequality always hold by Remark \ref{remIIA}.
Note that, for it construction, $x_J$ has no zero entries (in $J$). 
By the previous case, 
$x_J$ is an eigenvector of $B_J\circ x_Jx_J^*$ with eigenvalue $\IMBJ$. 
But clearly $B\circ xx^* $ has zeroes outside $J \times J$,  
so $x$ is an eigenvector of $B\circ xx^* $ iff $x_J$ is an eigenvector 
of $B_J\circ x_Jx_J^*$. Note that the eigenvalue of $x$ 
is always $\IMBJ$ \QED

\bigskip
\begin{teo}\label{elteo} Let $A \in P(n)$. Then  
$$I(2, A) = I(sp ,  {\bar A}\circ A)^{1/2}.
$$
\end{teo}
\dem 
Denote by $B = {\bar A}\circ A$. Given $y \in \Bbb C  ^n$ with $\|y\|=1$, 
we have that 
$$
\|A \circ yy^*\|_{_2}^2 = \sum_{i,j} |a_{ij}|^2 |y_i|^2 | y_j|^2 = 
\ \<  (B\circ yy^*) y ,  y\>  \ 
\le \|B\circ yy^*\|.
$$ 
Therefore $I(2, A)^2 \le I(sp , B)$. On the other hand, let $x$ 
be a unit vector with nonnegative entries such that $I(2,A)^2 = 
\|A \circ xx^*\|^2 $ and $ J= \{ i : x_i \ne 0 \}$.
Then, by Proposition \ref{el x},  $(B\circ xx^*)x = \IMBJ x$ and 
$$
I(2,A)^2 = \|A \circ xx^*\|^2 =  \ \<  (B\circ xx^*)x , x\>  \ = \IMBJ .
$$
But $x_J$ is a unit eigenvector of $B_J\circ x_Jx_J^*$ with strictly 
positive entries. So, by Lemma \ref{autovalor},  $B_J  (x_J \circ x_J) = 
\IMBJ (1, \dots , 1)^* $. Suppose that $I(2, A) \ne 0$. 
Then $\IMBJ \ne 0$ and,  by Theorem
\ref{i1=i2} and Remark \ref{remIIA}, 
$$
\IMBJ = I(sp, B_J)\ge I(sp, B) \ge I(2, A)^2 = \IMBJ .
$$
If $I(2,A) = 0$, then by Remark \ref{Rem1} some $A_{ii} = 0$. So also 
$I(sp, B) = 0$ by the same Remark \ \ \QED

\medskip
\begin{cor}\label{D2}
Let $A \in P(n)$. Then 
$$\begin {array}{rl}
I(2,A) &= \inf \ \{ \  (\sum _1^n D_{ii}^{-2})^{-1/2} : \ 0< D \ \, \mbox{ \rm 
is diagonal and } \, A \circ {\bar A} \le D^2 \ \}\\&\\
       & = \inf \ \{ \  I(2, D)\ : \ 0< D \ \, \mbox{ \rm 
is diagonal and } \, A \circ {\bar A} \le D^2 \ \}.\end{array}
$$
\end{cor}
\dem It is a direct consequence  of  Theorem \ref{elteo} and 
Proposition 3.2 of \cite {St}. See also Corollary \ref{losp} and
Remark \ref{infs} below \ \QED

\medskip
\begin{rem}\label{nopos} \rm
Theorem \ref{elteo} was formulated in order to get information about a matrix $A\in P(n)$
using the matrix $B = {\bar A} \circ A \in P(n)$. But it can also be interpreted in the 
converse way, i.e. to get information 
about a matrix $B \in P(n)$ with nonnegative entries using
the matrix $A= (B_{ij}^{1/2})$. Unfortunately it may certainly happen that 
$A \notin P(n)$ and one should check that $A \in P(n)$ in order to use 
the Theorem in this way. Nevertheless this restriction can be removed in the following
way: Given a selfadjoint (but not necessarily positive) matrix $A\in M_n$, let us 
still consider the index
$$
I(2, A) = \min \{\  \|A \circ xx^* \|_2 \ :  \ \|x \| = 1 \ \}
$$
defined by just acting on rank one projections. This definition is consistent for
positive matrices by Proposition \ref{uno}. 

A careful observation of the proofs of Proposition \ref{el x} and Theorem \ref{elteo}
shows that they remain true using this new index 
if we replace the condition ``$A \in P(n)$''   by 
``$A=A^*$ and $B= {\bar A} \circ A \in P(n)$''. Note that 
Lemma \ref{el z}, Lemma \ref{autovalor} and Theorem
\ref{i1=i2} are only applied  
to the positive matrix $B$ or its principal submatrices. 
The inequality $I(2,A) \le I(2, A_J)$ of 
Remark \ref{remIIA} (which is also used in the proofs) remains valid for this new index.
This observation is useful in order to avoid the
unpleasant condition $A = (B_{ij}^{1/2} ) \in P(n)$ in the following result.
\end{rem}

\medskip
\begin{cor}\label{Bj}
Let $B \in P(n)$ such that $b_{ij}\ge 0$ for all $i ,  j$. Then there exists a 
subset $J_0$ of $\{ 1, 2, \dots , n\}$ such that 
$ I(sp , B) = I(sp , B_{J_0}) = I(B_{J_0}) $. Therefore 
$$
I(sp , B) = \min \{ \ I(sp , B_J)\ : \ 
I(sp , B_J) = \IMBJ \ \}.
$$
If $A =(b_{ij}^{1/2})$  (which may be not positive), 
then $J_0$ can be also characteri\-zed as $J_0 = \{ i : x_i \ne  0\}$ for  
some unit vector
$x$ such that $I(2, A) = \|A\circ xx^*\|_2$. Also $I(sp , B)
= \|B\circ xx^*\|=  \ \< B y ,y\>  $, where $y = (|x_1|^2, \dots , |x_n|^2)^*$ .
\end{cor}

\noi
\dem  
Use Remark \ref{nopos} and the proof of Theorem \ref{elteo} \QED


\section{An Operator Inequality} In this section we compute the indexes of
a particulr type of matrices and, as an appplication, we get a new operator 
inequality, closely related to the inequality proved in \cite{CPR}, 
(see also \cite{ACS}, \cite{lo}).

Let $x=( \la_1 , \dots , \la _n )^*\in \Bbb R_+^n$,  
$S=\{ \la_1 , \dots , \la _n \}$ and  
$$
\Lambda = \Lambda_x = \Big{(}\la _i \la _j + {1\over \la _i  \la _j } \ \Big{)}_{ij}
\in P(n) .
$$

\medskip
\begin{num}\label{Lambda}\rm  {\bf Computation of $I(\La )$}

\medskip

\noi {\bf 1.} \ \ 
If all $\la _i$ are equal, then $\Lambda = (\la_1 ^2 +\la_1^{-2} ) \ P$ and $I(\La ) = 
\la_1 ^2 +\la_1^{-2} $.

\noi {\bf 2.} 
If $\# S >1$, then the image of $\Lambda$ is generated by
the vectors $x = (\la_1 , \dots , \la _n )^*$ and $y = (\la_1 \inv, \dots , \la _n\inv )^*$,
since $\Lambda = xx^* + yy^*$ and the matrix 
$$\Lambda_0 = \bm {cc} \la ^2 +\la^{-2} & \la \mu + \la \inv \mu \inv 
\\  \la \mu + \la \inv \mu \inv & \mu ^2 +\mu^{-2} \em
$$ 
is invertible if $\la \ne \mu$, so  rk $\Lambda  = 2$.

\noi {\bf 3.}  
If $\# S=2$, say $S= \{ \la , \mu \}$, 
then $p = a x + b y$, with $a = (\la + \mu)\inv $ and $b = \la \mu (\la + \mu)\inv $.
If a vector $z$ verifies that $\Lambda z = p$, then 
$$p = \Lambda z = (xx^* + yy^* )z = \ \< z,x\> x+\< z,y\> y. 
$$
Therefore
$$
I(\Lambda ) = \ \< z,p\> \inv = (\< z,x\> ^2 +\< z,y\> ^2 )\inv = 
{(\la + \mu)^2 \over 1+\la^2 \mu^2} = I(\Lambda_0) ,$$
where the last equality is shown in Remark 4.3 of \cite{St}.

\noi {\bf 4.}  
If $\# S >2$, it is easy to see that 
$p$ can not live in the subspace generated by $x$ and $y$. Then 
$I(\Lambda)$ 
must be zero  by Proposition \ref{proIIA}.
\medskip

\noi
Note that $I(\Lambda) \ne 0 $ if and only $\# S \le 2$.
\end{num}

\medskip
\begin{num}\label{Lambda2} \rm {\bf Computation of $I(sp, \La )$}
\medskip

\noi
We shall compute $I(sp, \La )$ using Corollary \ref{Bj} and therefore
use the principal minors of $\La$, which are
matrices of the same type. Let $J \subset \{ 1, 2, \dots , n\}$, 
$S_J = \{ \la_j : j \in J\}$ and $x_J$ is 
the induced vector. Then $\La_J = \La_{x_J}$ and so $I(sp, \La _J) \ne 0$. 
Suppose that $I(sp, \La _J) = I( \La _J)$. Then $\# S_J \le 2$ by \ref{Lambda}. 
If $\# S_J = 2$, let $i_1,i_2 \in J$ such that $\la_{i_1} \ne \la _{i_2}$. By
Theorem \ref{i1=i2} there exists a vector $0\le y \in \Bbb R^J$ 
such that  $\La_J \ y = p_J $. Let $z_1 = \sum \{ y_k : k \in J $ and $ 
\la_k = \la_{i_1} \}\ge0$ and $z_2 = \sum\{ y_j : j \in J $ and 
$ \la_j = \la_{i_2}\}\ge 0$. Easy computations show that 
$\La _{\{i_1,i_2\}} (z_1, z_2)^* = 
(1,1)^*$. Then, by Theorem \ref{i1=i2} and \ref{Lambda}, 
$$
I(sp, \La _J) = I(\La_J ) = 
{(\la_{i_1} + \la_{i_2})^2 \over 1+\la_{i_1}^2 \la_{i_2}^2}
= I( \La _{\{i_1,i_2\}})= I(sp, \La _{\{i_1,i_2\}}).
$$
Therefore, in order to compute $I(sp, \La )$ using Corollary \ref{Bj}, we only
have to consider the diagonal entries of $\La$ and some of the principal minors of size 
$2\times 2$. If $\la_i \ne \la_j$, by equation (\ref{bmin}),   
$$
I(sp, \La _{\{i,j\}}) = I( \La _{\{i,j\}}) \quad \Leftrightarrow \quad
\la_i \la_j +{1 \over \la_i \la_j}  \le \min \{\la_i^2+ {1 \over \la_i^2} \ , \ 
\la_j^2 + {1\over\la_j^2}\}.
$$
Easy computations show that, if $\la_i < \la_j$,  this condition 
is equivalent to 
\begin{equation}\label{lamu}
\la_i^2 \le {1\over \la_i \la_j} \le \la_j^2 .
\end{equation}
In particular, this implies that $\la_i < 1 < \la _j$. 
So, by Corollary \ref{Bj}, 
\begin{equation}\label{ILam}
I(sp, \La ) = \min \{M_1, M_2\}
\end{equation}
where $M_1 = \min_i \la_i ^2 + \la_i ^{-2}  = \min_i \La_{ii} $ \  and 
$$M_2 = 
\inf \Big{\{} {(\la_{i} + \la_{j})^2 \over 1+\la_{i}^2 \la_{j}^2} \ : \ 
\la_i < 1 < \la_j \ \hbox{ and } \ 
\la_i^2 \le {1\over \la_i \la_j } \le \la_j ^2 \Big{\}} .
$$
For example, if all $\la_i \ge 1 $ (or all $\la_i \le 1 $), then by equation (\ref{lamu}) 
$I(sp, \La )= M_1 = \min_i \la_i ^2 + \la_i ^{-2} $. 
On the other hand, 
if $\la \ne 1$ and $x = (\la , \la \inv )^*$, then 
$$
I(sp , \La _x ) = M_2 = {\la ^2 + \la ^{-2} \over 2}+1  <  M_1 = \la ^2 + \la ^{-2}.
$$
\end{num}

\medskip
\begin{pro}
Let ${\cal H}$ be a Hilbert space and $S$ a bounded selfadjoint invertible operator 
on ${\cal H}$. Let $M(S)$ be the best constant such that 
$$
\|STS + S\inv T S\inv \| \ge M(S) \|T\| 
\quad \hbox{ for all }    \ 0\le T \in L({\cal H}) 
$$
Then  $M(S) = \min \{ M_1 (S) , M_2 (S)\}$, where
$$M_1(S) = \min_{\la \in \ \sigma (S)} \la ^2 + \la ^{-2}
\quad \quad and  
$$
$$ M_2(S)  =  \inf   \Big{\{} {(|\la | + |\mu |)^2 \over 1+\la^2 \mu^2} :
\la, \mu \in \sigma (S) , |\la | < |\mu | \hbox{ and }  
\la^2 \le {1\over |\la \mu |} \le \mu ^2 \Big{\}} , 
$$
In particular, 
if $\|S\|\le 1 $ (or $\|S\inv \|\le 1$), then 
$$M(S) = \|S\|^{2} + \|S \|^{-2} \quad  
(resp.  \ \ \|S\inv \|^{2} + \|S\inv \|^{-2}).$$
\end{pro}

\noi\dem
We shall use the same steps as in \cite {CPR}. 
By taking the polar decomposition of $S$,  we can suppose that $S> 0$, since the unitary 
part of $S$ is also the unitary part of $S\inv$, 
commutes with $S$ and $S\inv$ and doesn't change norms. Note 
that we have to change $\sigma (S)$ by $\sigma(|S|) = 
\{|\la |: \la \in \sigma (S)\}$.

By the spectral theorem, we can suppose that $\sigma (S)$ is finite, 
since $S$ can be approximated in norm by operators $S_n$ such that each $\sigma (S_n)$ 
is a finite  subset of $\sigma (S)$,  $\sigma (S_n) \subset \sigma (S_{n+1})$
for all $n \in \zN$ and $\cup _n \sigma (S_n) $ is dense in $\sigma (S)$. 
So $M(S_n)$ (and $M_i(S_n), \ i = 1,2$) converges to $M(S)$ (resp. $M_i(S), \ i = 1,2$).

We can suppose also that $\dim {\cal H} < \infty$, by choosing an appropiate net of finite 
rank projections $\{P_F\}_{F \in {\cal F}}$ which converges strongly to the identity
and replacing $S, \ T$ by $P_F S P_F, \ P_F T P_F$. Indeed, the net may be choosen 
in such a way that $SP_F = P_F S$ and 
$\sigma (P_F S P_F ) = \sigma (S)$ for all $F\in {\cal F}$. Note that
for every $A\in L(H)$, $\|P_F A P_F\|$ converges to $\|A\|$.

Finally, we can suppose that $S$ is 
diagonal by a unitary change of basis in $\Bbb C ^n$. 
In this case, if $\la_1, \dots , \la _n$ are the eigenvalues of $S$ (with multiplicity) 
and $x = (\la_1, \dots , \la _n)^*$, then
$$
STS + S\inv T S\inv = \La _x \circ T .
$$
Note that all our reductions (unitary equivalences and
compressions) doesn't change the fact that $0\le T$.  
Now the statement follows from  formula (\ref{ILam}). If $\|S\| \le 1$ then 
$M(S) = M_1(S)$, since $M_2(S)$ is the infimum of the empty set. Clearly $M_1(S)$ 
is attained at the element $\la \in \sigma (S)$ such that $|\la | = \|S\|$ \quad \QED

\section{General unitary invariant norms}

\noi
Let $N$ be an unitary invariant norm in $M_n$ and let 
$\Phi$ be the symmetric gauge function on $\Bbb R^n$ associated with $N$
(see, for example, Chapter IV of \cite{Bh}). 
Our formulae are closely related to the mean 
$$
M_\Phi : ({\Bbb R} ^+ ) ^{n} \to \Bbb R ^+ , \quad \hbox{ given by } \quad  
M_\Phi (d_1, \dots ,d_n) = \Phi ' (d_1 \inv , \dots ,d_n\inv)\inv ,
$$
where $\Phi '$ is the dual norm of $\Phi$. Of particular interest are
those means induced by the Schatten $p$-norms, i.e 
$M_p: ({\Bbb R} ^+ ) ^{n} \to \Bbb R ^+$, $1< p \le \infty $ given by
$$
M_p(d_1, \dots ,d_n) = (\sum_{i=1}^n d_i^{-q})^{-1/q} 
=\|(d_1 \inv, \dots ,d_n \inv)\|_q\inv
 , \quad  (d_1, \dots , d_n) \in (\Bbb R ^+ )^n,
$$
where $q$ is the conjugate number of $p$. Coherently, we define $M_1 $ by
$$
M_1 (d_1, \dots ,d_n) = \min _{1\le i \le n} d_i .
$$
Note that ${1 \over n} M_\infty $ is usually known as the \it harmonic mean. \rm
\bigskip

In the following Remark we state several elementary properties of
the index $I(N, \cdot\ ) $ which hold for every unitary invariant norm: 

\medskip
\begin{rem}\label{Rem1} \rm Let $A \in P(n)$ and $N$ an  
unitarily invariant norm. 
Then
\ben
\item If we replace $N$ by $m \ N$ for some $m>0$, then the associated \ha \ index 
does not change. Therefore, from now on we shall assume that the norms are 
normalized, i.e. $N(E_{11}) = 1$.
\item $I(N, A) \le \min _{1\le i \le n} A_{ii}$, which can be seen  just
using the matrices $E_{ii}$.
\item Note that $\| \cdot \|_{sp}\le N(\cdot ) \le  \| \cdot \|_1 \le n \| \cdot \|_{sp}$. 
Then, if $A $ has no zero diagonal entries,  
$$
(\sum _1^n A_{ii}^{-1})^{-1} \le n \ I(N, A),
$$
since $(\sum _1^n A_{ii}^{-1})^{-1} \le I(sp, A)$, by  Corollary 3.7 of \cite{St}.

\item A consequence of the last inequalities is that 
$$I(N, A) = 0 \quad \Leftrightarrow \quad \hbox{ some } A_{ii} = 0 .
$$
\item If $A \le B$, then $I(N, A) \le I(N, B)$, since 
the same inequalities hold for the singular values of $A$ and $B$.
\item If $J \inc \{ 1, \ 2, \dots, \ n\}$ and $A_J$ is the principal 
submatrix of $A$ associated to $J$. Then $I(N, A)\le I(N, A_J)$. Indeed, 
the minimum which defines $I(N, A_J)$ is taken over less matrices
than the minimum which defines $I(N, A)$.
\een
\end{rem}

\medskip
\begin{pro} \label{ran1}
Let $A \in P(n)$ a rank one matrix. Then, for every
unitary invariant norm $N$,  
$$
I(N,A) =  M_1(A_{11} , \dots , A_{nn} ) =
\min_{1\le i\le n} A_{ii} ,
$$
\end{pro}
\dem Let $x \in \Bbb C  ^n $ such that $A = xx^*$. Let 
$$
m = \min_{1\le i \le n} |x_i|^2 = (\min_{1\le i\le n} A_{ii}) .
$$
We need to show that for every $B \in P(n)$, 
\begin{equation}\label{m2}
N(xx^*\circ B) \ge  m \ N(B),
\end{equation}
since this would imply that $I(N, A) \ge m$ and the other 
inequality always holds by 1. of Remark \ref{Rem1}. 
Clearly we can suppose that $m \ne 0$, 
so $x_i \ne 0$ $1\le i \le n$. 
Let $y = ( x_1\inv , \dots, x_n\inv)^t$. Then   
$B = yy^* \circ  (xx^* \circ B)$. 
Using a basic inequlity for the singular values (namely $s_i (\cdot)$) of a \ha 
product (Theorem 1 of \cite {AHJ}), we get that for $1\le k \le n$, 
$$
\begin{array}{rl}
\sum_1^k s_i(B)& = \sum_1^k s_i(yy^* \circ  (xx^* \circ B)) \\&\\
& \le \sum_1^k |{\tilde y}_i|^2 s_i(xx^* \circ B) \\&\\
& \le m\inv \sum_1^k  s_i(xx^* \circ B) , \end{array}
$$
where  ${\tilde y}$ is the vector $y$ rearranged in such a way that 
$|{\tilde y}_1| \ge |{\tilde y}_2| \ge \dots \ge |{\tilde y}_n|$
(which are the ``Euclidean norms'' of the ``columns'' of the $1\times n$ matrix $y^*$).
But this implies that equation (\ref{m2}) 
is true for the Ky-Fan norms $\| \cdot\|_{(k)}$ 
and therefore for every unitary invariant norm (see Theorem IV.2.2 of \cite{Bh})
\QED

\medskip
\begin{pro}\label{lasN}
For every unitary invariant norm $N$ and every 
diagonal matrix $D >0$, 
$$
I(N, D) = N'(D\inv) \inv , 
$$
where $N'$ is the dual norm of $N$. 
\end{pro}

\noi \dem 
Denote by $d_i = D_{ii}$, $d= (d_1, \dots , d_n)^*$ and 
$d\inv= (d_1\inv, \dots , d_n\inv)^*$.
Let $\Phi$ be the symmetric gauge function on $\Bbb R^n$ associated with $N$. 
Let $\Phi'$ be the dual  norm of $\Phi$ on $\Bbb R^n$. Then $\Phi'$  
corresponds to the dual norm $N'$ in $M_n$ (see IV.2.11 of \cite{Bh})
and, for any $A\in P(n)$, 
$$
\begin{array}{rl}
N( A ) \le \|A\|_1& =  \sum_1^n A_{ii} = \sum_1^n d_i \inv d_i A_{ii}  \\&\\
           & \le \Phi'( d\inv ) \Phi( d_i A_{ii}) \\ & \\
		   & = N'(D\inv)N(D\circ A ).
		   \end{array}
$$
Therefore $I(N, D) \ge N'(D\inv) \inv $. On the other hand, 
let $y \in \Bbb R ^n$  such that 
$\Phi  (y ) =1  $ and $ \<  d\inv  , y \>  \ = \Phi '(d\inv) $. 
Since $\Phi$ is a symmetric gauge function, we can suppose 
that all $y_i \ge 0 $. Let 
$x = ( d_1^{-{1\over2}} (y_{_1})^{{1\over2}}, 
\dots , d_n^{-{1\over2}} (y_{_n})^{{1\over2}})^*$. Then 
$$
 N(xx^*) = \|xx^* \|_1 = \hbox{tr } (xx^* ) =\ \<  d\inv  , y \>  \ = \Phi '(d\inv)
$$
and
$$
N( D \circ xx^* ) = \Phi (y ) = 1.
$$
Therefore 
$$
I(N, D) \le {N(D\circ xx^*) \over N(xx^*)} = \Phi '(d\inv)\inv =
N'(D\inv )\inv \quad \QED
$$

\medskip
\begin{rem}\label{cotinf} \rm
Given $A\in P(n)$ and an unitary invariant norm $N$, 
the lower bound for $I(N, A)$ given in 2. of Remark \ref{Rem1} can
be improved in the following manner: denote $D(A) = A \circ I$ and suppose that 
$D(A)$ is invertible. Then
$$
I(N, A) \ge I(N, D(A)) = N'(D(A) \inv ) \inv 
= \Phi' (A_{11}\inv , \dots , A_{nn}\inv ) \inv > 0,
$$
with $N'$ and $\Phi'$ as in the proof of Proposition \ref{lasN}. 
Indeed, it is well known that for every matrix $C$, \ 
$N(I\circ C) \le N(C)$ (for example, because $I\circ C$ is a convex combination
of matrices unitarily equivalent to $C$). Therefore for every $B\in P(n)$, 
$$
N(A\circ B) \ge N(I\circ A \circ B) = N(D(A)\circ B) \ge I(N, D(A) )\  N(B)
$$
\end{rem}

\medskip
\begin{cor}\label{losp}
Let $0< D $ be a diagonal matrix. Then, 
for $1 < p \le \infty$ and $\frac1p + \frac1q = 1$, 
we have that 
$$
I(p, D) = M_p(D_{11} , \dots , D_{nn} ) = (\sum _1^n D_{ii}^{-q})^{-1/q}.
$$
\end{cor}


\medskip
\begin{rem}\label{infs} \rm In the case of the spectral norm, the index of the
diagonal matrices determine the index of all positive matrices. Indeed, 
it is shown in Proposition 3.2  of \cite{St} that 
\begin{equation}\label{inf}
\begin{array}{rl}
I(sp, A) &
= \inf  \ \{ \ I(sp, D) : A\le D \ \hbox{ and $D$ is diagonal }\}
\end{array}
\end{equation}
since the condition $\bm{cc} D&A\\A&D\em \ge 0$ which appears in that Proposition
is clearly equivalent to $A \le D$. One could be tempted to conjecture that
similar formulae hold for other norms. Unfortunately this assertion fails, 
at least in this simple formulation. For example, Corollary \ref{D2} says 
that, using the Frobenius norm $\| \cdot\|_2$, for every $A \in P(n)$, 
\begin{equation}\label{inf2}
\begin {array}{rl}
I(2,A) &= \inf \{ (\sum _1^n D_{ii}^{-2})^{-1/2} : D \mbox{ \rm 
is diagonal and }  A \circ {\bar A} \le D^2 \}\\&\\
       & = \inf \ \{ \  I(2, D)\ : \ 
	   D \ \, \mbox{ \rm 
is diagonal and } \, A \circ {\bar A} \le D^2 \ \}.\end{array}
\end{equation}
Note that the condition $A \circ {\bar A} \le D^2$ is strictly less restrictive
that $A\le D$ (the reverse implication follows from Schur Theorem). 
It can be easily seen with a computer that formula (\ref{inf}) 
does not hold for the Frobenius norm. Nevertheless,
equation (\ref{inf2}) allows one to compute the 2-index for every positive matrix
using only diagonal matrices. 
\end{rem}

\end{document}